\documentclass{amsart}

\usepackage{pstricks, pst-tree, pst-node} 
\usepackage[latin1]{inputenc}
\usepackage{amssymb}
\usepackage{xfrac}     
\usepackage{lmodern}
\usepackage{stmaryrd}
\usepackage{graphicx}
\usepackage{amsmath}  
\usepackage{amsthm}
\usepackage{wrapfig} 
\usepackage{color}  
\usepackage{microtype} 
\usepackage{setspace}
\usepackage{amsfonts}
\usepackage{wasysym}
\usepackage{cancel}
\usepackage{verbatim}
\usepackage{dsfont}
\usepackage{paralist}
\usepackage{mathtools}

\DeclareMathOperator{\spann}{span}

\DeclareMathOperator{\FF}{\mathbb{F}}

\DeclareMathOperator{\ZZ}{\mathbb{Z}}

\DeclareMathOperator{\NN}{\mathbb{N}}

\DeclareMathOperator{\w}{\wedge}
\DeclareMathOperator{\de}{d\!}
\DeclareMathOperator{\Ker}{Ker}
\DeclareMathOperator{\Imm}{Im}
\DeclareMathOperator{\pB}{\mathcal B}
\DeclareMathOperator{\coker}{coker}

\DeclareMathOperator{\pdeg}{\textit{p}\hspace{1.5pt}-deg}

\DeclareMathOperator{\ann}{Ann}

\newcommand{\be}{\hfill$\square$}

\newtheoremstyle{plain2}
  {10pt}   % ABOVESPACE
  {10pt}   % BELOWSPACE
  {\itshape}  % BODYFONT
  {0pt}       % INDENT (empty value is the same as 0pt)
  {\bfseries} % HEADFONT
  {}         % HEADPUNCT
  {5pt plus 1pt minus 1pt} % HEADSPACE
  {}          % CUSTOM-HEAD-SPEC
  
 \newtheoremstyle{beweis}
  {10pt}   % ABOVESPACE
  {10pt}   % BELOWSPACE
  {\normalfont}  % BODYFONT
  {0pt}       % INDENT (empty value is the same as 0pt)
  {\bfseries} % HEADFONT
  {:}         % HEADPUNCT
  {5pt plus 1pt minus 1pt} % HEADSPACE
  {}          % CUSTOM-HEAD-SPEC

\newtheoremstyle{definition2}
  {10pt}   % ABOVESPACE
  {10pt}   % BELOWSPACE
  {\normalfont}  % BODYFONT
  {0pt}       % INDENT (empty value is the same as 0pt)
  {\bfseries} % HEADFONT
  {}         % HEADPUNCT
  {5pt plus 1pt minus 1pt} % HEADSPACE
  {}          % CUSTOM-HEAD-SPEC

\theoremstyle{plain2}
\newtheorem{satz}{Satz}[section]
\newtheorem{lem}[satz]{Lemma}
\newtheorem{pro}[satz]{Proposition}
\newtheorem{theo}[satz]{Theorem}
\newtheorem{coro}[satz]{Corollary}

\theoremstyle{definition2}
\newtheorem{defi}[satz]{Definition}
\newtheorem{rem}[satz]{Remark}

\newtheorem{bsp}[satz]{Example}

\theoremstyle{beweis}
\newtheorem*{prf}{Proof}

\begin{document}

\title[Annihilators of differential forms]{Annihilators of differential forms over fields of characteristic $p$}

\author{Marco Sobiech}
\address{Fakult\"at f\"ur Mathematik, Technische Universit\"at Dortmund, D-44221 Dortmund, Germany}
\email{marco.sobiech@tu-dortmund.de}
\date{\today}

\begin{abstract}
Let $F$ be a field of characteristic $p$ and let $\Omega^n(F)$ be the $F$-vector space of $n$-differential forms. 
In this work, we will study the annihilator of differential forms, 
give specific descriptions for special cases and show a connection between these annihilators and 
the kernels of the restriction map $\Omega^n(F) \to \Omega^n(E)$ for purely inseparable field extensions $E/F$.
\end{abstract}



\maketitle

\begin{section}{Introduction}

For a field $F$ of characteristic $p>0$, the space of $n$-fold differential forms $\Omega^n(F)$ is defined as the $n$-fold exterior power of the $F$-vector 
space $\spann_F(\de F)$. The Artin-Schreier map $\wp$ naturally extents to $\Omega^n(F)$, so in addition to the space of $n$-differential forms, we also 
get two more spaces $ \nu_n(F):=\ker\wp$ and $H^{n+1}_p(F):=\coker \wp$. Since the structures of these three spaces 
heavily depend on the field $F$, it is a natural question to ask how these spaces behave under field extensions $K/F$. 
In particular we are interested in a characterization of forms defined over $F$ that 
become the zero form viewed as an element over the field $K$, i.e. we are interested in the kernel of the restriction map induced by $F \to K$ which we will 
denote by $\Omega^n(K/F)$, $\nu_n(K/F)$ and $H^{n+1}_p(K/F)$. 
In this note we will focus on algebraic extensions of $F$. Using the theory of $p$-independence, it is well known that $\Omega^n(K/F)$ and $\nu_n(K/F)$ are both trivial 
for all separable extensions $K/F$. So the next type of field extensions to study which come to mind are 
purely inseparable extensions. In this case the kernel $H_p^{n+1}(K/F)$ was completely described in \cite{q55}. For the kernels $\Omega^n(K/F)$ and $\nu_n(K/F)$ 
not so much is known. For modular purely inseparable extensions given by $K=F\left( \sqrt[p^{m_1}]{b_1},\ldots,\sqrt[p^{m_r}]{b_r}\right)$ (i.e. the set $\{ b_1,\ldots,b_r\}$ is $p$-independent), in \cite[Theorem 4.1]{q55} was shown  
$$\Omega^n(K/F)=\sum_{i=1}^r \de b_i \w \Omega^{n-1}(F)$$
and if we additionally assume $F^{p-1}=F$ so that we can use the Lemma of Kato 
%stated here as Lemma \ref{23} 
we get 
$$\nu_n(K/F)= \left[ \frac{\de x}{x} \mid x \in F^p(b_1,\ldots, b_r)^* \right] \w \nu_{n-1}(F)$$
(see \cite[Theorem 5.7]{q55} for $p=2$, but the proof can easily be translated to arbitrary positive characteristic). The case of non modular purely inseparable extensions however seems to be more difficult. In Theorem \ref{54} we will provide insight into the 
kernels $\Omega^n(K/F)$ and $\nu_n(K/F)$ for a special case of non modular purely inseparable extensions which in particular gives us a complete characterization for 
$2$-fold purely inseparable extensions, i.e. for purely inseparable extensions generated by two elements. Our main tool in the study of these kernels are so called 
annihilators, which we will introduce in sections 2 and 3 to proof Theorem \ref{43}, so that we can translate kernels into annihilators. 
It turns out that we can describe all 
kernels $\Omega^n(K/F)$ and $\nu_n(K/F)$ which are known at the time of publishing as annihilators of specific subsets of $\Omega^n(F)$ which contain the crucial information of the algebraic field extension $K/F$.
%It turns out that by using these annihilators, every kernel $\Omega^n(K/F)$ and $\nu_n(K/F)$ 
%known so far may be written 
%as an annihilator of a specific subset of $\Omega^n(F)$ which contains the crucial information of the field extension.

\end{section}

\begin{section}{A short introduction to differential forms}

We refer to \cite{q22} or \cite{q1} (for the case $p=2$) for any undefined terminology 
or any basic facts about
differential forms that we do not mention explicitly in this work.
%We will assume $0\in \NN$.

Throughout this whole paper let $F$ denote a field of characteristic $p>0$, if not stated otherwise.
We start by recalling some basic facts about the theory of differential forms. 

The space $\Omega^1(F)$ of absolute $1$-differential forms over $F$ is defined to be the $F$-vector 
space generated by the symbols $\de a$ with $a \in F$ which are subject to the relations 
$\de\, (a+b)=\de a + \de b$ and $\de\,(ab)=a\de b + b \de a$ for all $a,b \in F$. 
With this we have $\de a^p=0$ for all $a\in F$ and $\de : F \to \Omega^1(F), a \mapsto \de a$ 
is an $F^p$-derivation. The space of $n$-fold differential forms, or $n$-differentials for 
short, is defined as the $n$-fold exterior power $\Omega^n(F):= \bigwedge^n \Omega^1(F)$ 
with $n \in \NN$. Therefore $\Omega^n(F)$ is an $F$-vector space generated by the 
elementary wedge products $\de a_1 \w \ldots \w \de a_n$ with $a_1,\ldots,a_n\in F$. With this, the
 map $\de$ can be extended to an $F^p$ linear map, defined on additive generators by 
 $$\de: \Omega^n(F) \to \Omega^{n+1}(F)\ ,\  x \de a_1 \w \ldots \w \de a_n \mapsto 
 \de x \w \de a_1 \w \ldots \w \de a_n.$$ 
 Forms contained in the image $\de \Omega^n(F)$
 of the operator $\de$ are called exact differential forms. For completeness of later 
 results, we further define $\Omega^0(F):=F$ 
and $\Omega^n(F):=\{0\}$ for $n<0$. The Artin-Schreier map can be extended 
to the space $\Omega^n(F)$ by 
\begin{align*}
\wp : \Omega^n(F) &\to \Omega^n(F) / \de \Omega^{n-1}(F)\ , \\
x \frac {\de a_1}{a_1}\w \ldots \w \frac{\de a_n}{a_n} &\mapsto
(x^p-x) \frac {\de a_1}{a_1}\w \ldots \w \frac{\de a_n}{a_n} \mod \de \Omega^{n-1}(F)
\end{align*}
with $\Ker (\wp)=:\nu_n(F)$ and Kato showed in \cite{q4}, that $\nu_n(F)$ is 
additively generated by the so called logarithmic differential 
forms $\frac {\de a_1}{a_1}\w \ldots \w \frac{\de a_n}{a_n}$ with $a_1,\ldots,a_n \in F^*$.

Now recall that a subset $\mathcal A\subset F$ is called $p$-independent (over $F$), 
if for every finite subset $\{a_1,\ldots,a_k\}\subset \mathcal A$, we 
have $[F^p(a_1,\ldots,a_k):F^p]=p^k$. Additionally $\mathcal A$ is called a 
$p$-basis of $F$, if $\mathcal A$ is $p$-independent with $F^p(\mathcal A)=F$ and we 
define the $p$-degree of a set $\mathcal C \subset F$ as $\pdeg(\mathcal C) := \log_p([F^p(\mathcal C) : F^p])$ if $[F^p(\mathcal C) : F^p]$ is finite and 
$\infty$ otherwise . With this we 
get the following lemma which is folklore by now, but we will state it here for easy reference.

	\begin{lem}\label{21}
	\begin{enumerate}[(a)]
		\item Let $a_1,\ldots,a_n\in F$. Then the following statements are equivalent
			\begin{enumerate}[(i)]
				\item $\{a_1,\ldots,a_n\}$ is $p$-independent.
				\item $\de a_1,\ldots,\de a_n$ are $F$-linearly independent.
				\item $\de a_1 \w \ldots \w \de a_n \neq 0$ in $\Omega^n(F)$.
			\end{enumerate}
		\item Let (a) be true for $a_1,\ldots,a_n\in F$ and let $b_1,\ldots,b_n\in F$. then the following statements are equivalent
		\begin{enumerate}[(i)]
			\item $F^p(a_1,\ldots,a_n)=F^p(b_1,\ldots,b_n)$.
			\item $\spann_F(\de a_1,\ldots, \de a_n)=\spann_F(\de b_1,\ldots,\de b_n)$.
			\item $\de a_1\w\ldots\w\de a_n= x \de b_1 \w \ldots \w \de b_n$ for some $x \in F^*$.
		\end{enumerate}
		
		\item $\de a \in \spann_F(\de b_1,\ldots,\de b_n)$ if and only if  
$a \in F^p(b_1,\ldots,b_n)$.
	\end{enumerate}
	\end{lem}

As an easy consequence from Lemma \ref{21} we get $\Omega^n(F)=\{0\}$ if $n>|\pB|$ for a $p$-basis $\pB$ of $F$. Using basic linear algebra, we can also 
deduce from Lemma \ref{21} that every $p$-independent set in $F$ can be extended to a $p$-basis of $F$. This seems to be a trivial observation but is rather crucial 
for many of the upcoming proofs.

For two differential forms $\omega\in \Omega^n(F)$ and $u\in \Omega^r(F)$ we say that 
$\omega$ divides $u$, for short $\omega \mid u$, if there exists a $v \in \Omega^{r-n}(F)$ 
with $u= \omega \w v$. Now since every $p$-independent subset can be extended to a $p$-basis of $F$, it is evident to see the following

	\begin{lem}\label{22}%{\cite[Prop. 1.16 for $p=2$]{q1}}
Let $\{b_1,\ldots,b_r\}\subset F$ be $p$-independent and $\omega\in \Omega^n(F)$. Then 
\begin{align*}
\de b_i \mid \omega\ \text{ for } i=1,\ldots,r\ \iff \de b_1 \w \ldots \w \de b_r \mid \omega
\end{align*}
	\end{lem}

Now let $\pB=\{ b_i \mid i\in I\}$ be a $p$-basis of $F$. We may assume that $I$ is well ordered and 
set $\Sigma_n:=\{ \sigma:\{1,\ldots,n\} \to I \mid \sigma(i)<\sigma(j)$ for $i<j \}$. We further transfer the ordering of $I$ to $\mathcal B$ by setting $b_i<b_j$ if $i<j$ in $I$. 
Then it 
is well known that the differential basis 
\begin{equation*}
\bigwedge\nolimits^n_{\pB}:=\{ \de b_{\sigma(1)}\w \ldots \w \de b_{\sigma(n)} \mid \sigma \in \Sigma_n \}
\end{equation*} 
%resp. 
%$\{\frac {\de b_{\sigma(1)}}{b_{\sigma(1)}}\w \ldots \w \frac{\de b_{\sigma(n)}}{b_{\sigma(n)}} \mid \sigma \in \Sigma_n\}$
is an $F$-basis of $\Omega^n(F)$. For $b \in \mathcal B$ and $\sigma \in \Sigma_n$ we will  shortly write $b\in \Imm (\sigma)$ if $b=b_i$ for some $i\in I$ and $i\in\Imm (\sigma)$.

Furthermore we can equip $\Sigma_n$ with the lexicographic ordering, i.e. for $\sigma,\tau \in \Sigma_n$ 
we have $\sigma < \tau$ if and only if $\sigma(j)< \tau(j)$ and $\sigma(i)=\tau(i)$ for 
some $j\in\{1,\ldots,n\}$ and all $i<j$. Setting $\de b_\sigma:=\de b_{\sigma(1)}\w \ldots \w \de b_{\sigma(n)}$ resp. 
$\frac{\de b_{\sigma}}{b_{\sigma}}:=\frac {\de b_{\sigma(1)}}{b_{\sigma(1)}}\w \ldots \w \frac{\de b_{\sigma(n)}}{b_{\sigma(n)}}$ 
for $\sigma \in \Sigma_n$, we get a filtration of $\Omega^n(F)$ by 
$$\Omega^n_{\delta}(F):= \spann_F\left( \de b_{\sigma} \mid \sigma \leq \delta \right)\ \ \text{ resp. } \ \ 
\Omega^n_{<\delta}(F):= \spann_F\left( \de b_{\sigma} \mid \sigma < \delta \right),$$
with $\Omega^n_{\delta}(F)\subsetneq \Omega^n_{\tau}(F)$ if and only if $\delta < \tau$. 
Using this notation, we set 
$\max_{\pB}(\omega):=\min\{ \delta \in \Sigma_n \mid \omega \in \Omega^n_{\delta}(F)\}$ 
and call it the maximal multiindex of $\omega$ with respect to $\pB$. 
Note that $\wp$ respects this filtration on $\Omega^n(F)$ after fixing 
a $p$-basis of $F$ resp. the corresponding 
differential basis.
Transferring the above filtration to the field $F$, we define
$$ F_{\ell}:= F^p( b_i \mid i\leq \ell)\ \ \text{ and } \ \ F_{<\ell}:= F^p( b_i \mid i< \ell)$$
for any $\ell \in I$.

\end{section}

\begin{section}{Annihilators in $\Omega^n(F)$}\label{sec2}

We will now start to analyse annihilators in the algebra of differential forms. 
First, we make the following definition.

\begin{defi}\label{31}
For $n,r \in \NN$ and a non-empty set $U \subset \Omega^r(F)$, we define 
\begin{align*}
\ann \Omega^n_F(U):&= \{ \omega\in \Omega^n(F)  \mid \omega\w u =0 \in \Omega^{n+r}(F) \text{ for all } u\in U  \} \\
\ann \nu^n_F(U):&= \{ \chi\in \nu_n(F)  \mid \chi\w u =0\in\Omega^{n+r}(F) \text{ for all } u\in U  \}
\end{align*}
and call it the $\Omega$-annihilator resp. the $\nu$-annihilator of $U$.
\end{defi}
Note $\ann\nu^n_F(U)=\ann\Omega^n_F(U) \cap \nu_n(F)$ and also note $\ann\Omega^n_F(U)=\Omega^n(F)$ if $|\pB|<n+r$ for a $p$-basis 
$\pB$ if $F$. 
Our goal will be to study these annihilators for specific types of $U$, which are 
given as follows. For non-empty sets $S_1,\ldots,S_r\subset F$, we define 
$\de S_1 \w \ldots \w \de S_r= \{ \de s_1 \w \ldots \w \de s_r \mid s_i\in S_i, i=1,\ldots,r\}$. 
For $S_i=\{s_i\}$, we will replace the slot $\de S_i$ by $\de s_i$ for short and for a single non-empty set $S\subset F$, we 
will write $\bigwedge^r \de S$ for the $r$-fold wedge product $\de S \w \ldots \w \de S $.

\begin{lem}\label{32}
For non-empty sets $S_1,\ldots,S_r \subset F$ and $U,V \subset \Omega^r(F)$, we have
\begin{enumerate}[(a)]
	\item $\left( \de S_1 \w \ldots \w \de S_r \right) \cup \{0\} = \de\,( S_1\setminus F^p) \w \ldots \w \de\, (S_r \setminus F^p)\cup \{0\}$ if $S_i\setminus F^p \neq \varnothing$.
	\item $\ann\Omega^n_F(\de S_1 \w \ldots \w \de S_r)= \ann\Omega^n_F(\de S_{\pi(1)} \w \ldots \w \de S_{\pi(r)})$ for all permutations $\pi$.
	\item For $U \subset V$ we have $\ann\Omega^n_F(V) \subset \ann\Omega^n_F(U)$.
	\item $\ann\Omega^n_F(U) = \ann\Omega^n_F(\spann_F (U))$.
	\item $\ann\Omega^n_F \left( \de S_1 \w \ldots \w \de S_r \right) = 
	\ann\Omega^n_F \left( \de \,(F^p(S_1)) \w \ldots \w \de \,(F^p(S_r)) \right)$.

\end{enumerate}
\end{lem}

\begin{prf}
The statements (a)-(c) are all consequences from the definitions above and (d) is proved using the multilinear property of the wedge 
product. The inclusion ($\supseteq$) in (e) directly follows from (c). To check the inverse inclusion ($\subseteq$) let $\{ a_{ij_i} \mid j_i \in J_i\}$ be a $p$-basis of $F^p(S_i)$ for $i\in\{1,\ldots,r\}$. Using the Basis Extension Theorem for $p$-independence, we may assume $a_{ij_i}\in S_i$ for all 
$j_i \in J_i$. Now let $\omega\in \ann\Omega^n_F \left( \de S_1 \w \ldots \w \de S_r \right)$ and $x_i \in F^p(S_i)$. 
Then we find a finite subset $T_i \subset J_i$ and $\lambda_{i t_i}\in F$ for $t_i\in T_i$ with 
$\de x_i = \sum_{t_i \in T_i} \lambda_{it_i} \de a_{it_i}$. Thus we have 
\begin{equation*}
\omega \w \de x_1 \w \ldots \w \de x_i = 
\sum_{t_1\in T_1}\ldots\sum_{t_r\in T_r} \lambda_{1t_1}\ldots \lambda_{rt_r} \omega\w \de a_{1t_1}\w \ldots \w \de a_{rt_r}=0
\end{equation*}
due to the selection $a_{ij_i}\in S_i$. This concludes the proof. \be
\end{prf}

As an easy consequence, we get

\begin{coro}\label{33}
For non-empty sets $S_1,\ldots,S_r \subset F$ with $\pdeg(S_i)=k_i\in \NN$ and a $p$-basis $\{a_{i1},\ldots,a_{ik_i}\}$ of 
$F^p(S_i)$, we have 
\begin{align*}
\ann\Omega^n_F \left( \de S_1 \w \ldots \w \de S_r \right)= \ann\Omega^n_F \left( \de \,\{a_{11},\ldots,a_{1k_1}\} \w \ldots \w \de\, \{a_{r1},\ldots,a_{rk_r}\} \right).
\end{align*}
In particular we have 
\begin{align*}
\ann\Omega^n_F \left( \de S_1 \w \ldots \w \de S_r \right)= 
\bigcap_{j_1=1}^{k_1}\ldots \bigcap_{j_r=1}^{k_r} \ann\Omega^n_F( \de a_{1j_1}\w\ldots\w\de a_{rj_r})
\end{align*}
\end{coro}

With this setup, we are now able to compute our first annihilators.

\begin{pro}\label{34}
For non-empty sets $S_1,\ldots,S_r \subset F$ with $\pdeg(S_i)=k_i\in \NN$ and $p$-bases $\{a_{i1},\ldots,a_{ik_i}\}$ of $F^p(S_i)$ for $i=1,\ldots,r$ assume 
$\pdeg (S_1 \cup \ldots \cup S_r)=k_1+ \ldots +k_r$. Then 
\begin{align*}
\ann\Omega^n_F \left( \de S_1 \w \ldots \w \de S_r \right)= 
\sum_{i=1}^r \de a_{i1}\w\ldots\w\de a_{ik_i} \w \Omega^{n-k_i}(F).
\end{align*}
In particular, we have
	\begin{align*}
	\ann\Omega^n_F&(\de S_1 \w \ldots \w \de S_r)=\sum_{i=1}^r \ann\Omega^n_F( \de S_i) \\
	&=\sum_{i=1}^r \bigcap_{j_i=1}^{k_i}\ann\Omega^n_F(\de a_{ij_i}) 
	=\bigcap_{j_1=1}^{k_1}\ldots \bigcap_{j_r=1}^{k_r}\sum_{i=1}^r \ann\Omega^n_F(\de a_{ij_i}).
	\end{align*}

\end{pro}

\begin{prf}
As above, we may assume $a_{i1},\ldots,a_{ik_i}\in S_i$ for all $i\in\{1,\ldots,r\}$. Since $a_{11},\ldots, a_{1k_1},\ldots,a_{r1},\ldots,a_{rk_r}$ are 
$p$-independent by assumption, we can extend these elements to a $p$-basis $\pB=\{b_i \mid i\in I\}$ of $F$. All we have to prove is the 
equality 
\begin{align*}
\ann\Omega^n_F \left( \de S_1 \w \ldots \w \de S_r \right)= 
\sum_{i=1}^r \de a_{i1}\w\ldots\w\de a_{ik_i} \w \Omega^{n-k_i}(F)
\end{align*}
in which the inclusion $(\supseteq)$ is rather easy to see and left to the reader. So let us now prove the inclusion $(\subseteq)$.

Let $\omega \in \Omega^n(F) \setminus \sum_{i=1}^r \de a_{i1}\w\ldots\w\de a_{ik_i} \w \Omega^{n-k_i}(F)$. Thus in the basis expansion of $\omega$ using the 
basis $\bigwedge^n_{\mathcal B}$ of $\Omega^n(F)$, we find a non zero summand $\lambda = x \de b_{j_1} \w \ldots \w \de b_{j_n}$ with $x \in F^*, b_{j_1}, \ldots, b_{j_n}\in \mathcal B$, for which 
we can choose $a_i \in \{a_{i1},\ldots,a_{ik_i}\}$ such that $\de a_i \nmid \lambda$ for all $i\in\{ 1,\ldots,n\}$. Since all the elements  $a_{11},\ldots, a_{1k_1},\ldots,a_{r1},\ldots,a_{rk_r}$ are contained in the $p$-basis of $F$, by Lemma \ref{22} we get 
$$ \de a_1 \w \ldots \w \de a_n \w \lambda \neq 0$$
and therefore 
$$ \de a_1 \w \ldots \w \de a_n \w \omega \neq 0.$$
From this we can easily conclude that $\omega \in \Omega^n(F) \setminus \ann\Omega^n_F \left( \de S_1 \w \ldots \w \de S_r \right)$.\be
\end{prf}

\begin{pro}\label{35}
Let $S\subset F$ be a non-empty set with $\pdeg (S)=k$ and let $\{a_1,\ldots,a_k\}$ be a $p$-basis of $F^p(S)$. 
For $r\in\{1,\ldots,k\}$ set $t=k-r+1$. Then 
\begin{align*}
	\ann\Omega^n_F\left(\bigwedge\nolimits^r \de S \right) =\sum_{\substack{\{i_1,\ldots,i_t\}\subset\{1,\ldots,k\} \\ i_1<\ldots <i_t  }} 
	\de a_{i_1}\w\ldots \w\de a_{i_t} \w\Omega^{n-t}(F)
\end{align*}
and for $r>k$ we have $	\ann\Omega^n_F\left(\bigwedge\nolimits^r \de S \right)=\Omega^n(F)$.
\end{pro}

\begin{prf}
Using Corollary \ref{33}, we may assume $S=\{a_1,\ldots,a_k\}$. The second statement is clearly true since  $\bigwedge\nolimits^r \de S=\{0\}$ for $r>k$ by Lemma \ref{21}. 

So assume $r\in\{1,\ldots,k\}$ and set $t=k-r+1$. We can extend $\{a_1,\ldots,a_k\}$ 
to a $p$-basis $\pB=\{b_i \mid i\in I\}$ of $F$. Again the inclusion $(\supseteq)$ is easy to check, so we will now prove the 
inclusion $(\subseteq)$ and leave the rest to the reader. To simplify notations, set 
$$\mathfrak A =\sum_{\substack{\{i_1,\ldots,i_t\}\subset\{1,\ldots,k\} \\ i_1<\ldots <i_t  }} 
	\de a_{i_1}\w\ldots \w\de a_{i_t} \w\Omega^{n-t}(F)$$
	and let $\omega \in \Omega^n(F) \setminus \mathfrak A$. Thus the basis expansion of $\omega$ using $\bigwedge^n_{\mathcal B}$ contains a non 
	zero summand $\lambda = x \de b_{s_1}\w \ldots \w \de b_{s_n}$ such that $\de a_{i_1}\w\ldots \w\de a_{i_t} \nmid \lambda$ for all 
	subsets $\{i_1, \ldots, i_t\}\subset \{1,\ldots,k\}$. In particular we have $\de a_1 \w \ldots \w \de a_t \nmid \lambda$, hence we find 
	$j_1 \in \{1,\ldots,t\}$ with $\de a_{j_1} \nmid \lambda$ by Lemma \ref{22}. Similarly we find $j_2 \in \{1,\ldots,t+1\}\setminus\{j_1\}$ with 
	$\de a_{j_2}\nmid \lambda$ and continue until we have found $j_1,j_2,\ldots ,j_r\in\{1,\ldots,k\}$ such that $\de a_{j_1}\nmid \lambda, \ldots ,\de a_{j_r}\nmid \lambda$.
	Again by Lemma \ref{22} we then get $\de a_{j_1} \w \ldots \w \de a_{j_r} \nmid \lambda$ from which we conclude $\de a_{j_1} \w \ldots \w \de a_{j_r} \w \lambda \neq 0$ 
	and finally $\de a_{j_1} \w \ldots \w \de a_{j_r} \w \omega \neq 0$. 
\be
\end{prf}

\begin{pro}\label{36}
Let $S_1,\ldots,S_{r+1}\subset F$ be non-empty sets with $\pdeg S_i=1$ for $i=1,\ldots,r$ and $\de S_1 \w \ldots \w\de S_r \w \de S_{r+1}\neq \{0\}$. 
Set $\pdeg(S_1 \cup \ldots \cup S_r \cup S_{r+1})=r+\ell$ with $\ell \in \ZZ$. Then $\ell \geq 1$ and for any $p$-independent  
$a_1,\ldots,a_r,e_1,\ldots, e_{\ell}\in F$ with $F^p(S_i)=F^p(a_i)$ for $i=1,\ldots,r$ 
and $F^p(a_1,\ldots, a_r)(S_{r+1})=F^p(a_1,\ldots,a_r)(e_1,\ldots,e_{\ell})$ we have  
\begin{align*}
	\ann\Omega^n_F( \de S_1 \w \ldots \w \de S_r\w \de S_{r+1} ) = \sum_{i=1}^r \de a_i \w \Omega^{n-1}(F) + \de e_1\w\ldots \w\de e_{\ell}\w\Omega^{n-\ell}(F).
\end{align*}
\end{pro}

\begin{prf}
As above we may assume $S_i=\{a_i\}$ for $i=1,\ldots,r$ 
and $e_1,\ldots,e_{\ell}\in S_{r+1}$. First $\{a_1,\ldots,a_r\}$ must 
be $p$-independent, because a $p$-dependence would imply 
$\de S_1 \w \ldots \w \de S_r \w \de S_{r+1}=\{0\}$. 
So $\ell \geq 0$ and $\ell =0$ would 
imply $F^p(S_{r+1})\subset F^p(a_1,\ldots,a_r)$ which would lead 
to $\de S_1 \w \ldots \w \de S_r \w \de S_{r+1}=\{0\}$ as well. Thus we 
have $\ell \geq 1$ and $\{a_1,\ldots,a_r,e_1,\ldots, e_{\ell}\}$ is 
a $p$-basis of $F^p(S_1,\ldots,S_{r+1})$. 

Let us start with the proof of $(\subseteq)$. Since $\{e_1,\ldots,e_{\ell}\}$ can be seen as part of a $p$-basis $\mathcal S$ 
of $F^p(S_{r+1})$, by Corollary \ref{33},  Lemma \ref{32}(c) and Proposition \ref{34} 
we get
\begin{align*}
\ann\Omega^n_F(\de S_1 \w \ldots \w \de S_r \w \de S_{r+1}) &= 
\ann\Omega^n_F(\de S_1 \w \ldots \w \de S_r \w \de \mathcal S) \\
 &\subset \ann\Omega^n_F(\de S_1 \w \ldots \w \de S_r \w \de\, \{e_1,\ldots,e_{\ell}\})\\
 &=\sum_{i=1}^r \de a_i \w \Omega^{n-1}(F) + \de e_1\w\ldots \w\de e_{\ell}\w\Omega^{n-\ell}(F).
\end{align*}

For the reverse inclusion $(\supseteq)$ set $T:=\{a_1,\ldots,a_r\}\cup S_{r+1}$. Then we obviously have
\begin{align*}
\ann\Omega^n_F( \de S_1 \w \ldots \w \de S_r\w \de S_{r+1} ) = \ann\Omega^n_F( \de a_1 \w \ldots \w \de a_r\w \de T )
\end{align*}
and by Corollary \ref{33} we may replace $T$ by $\{a_1,\ldots,a_r,e_1,\ldots,e_{\ell}\}$ 
in the equation above. So it now suffices to show that the 
forms $\de a_1,\ldots, \de a_r$ and $\de e_1 \w \ldots \w \de e_{\ell}$ annihilate 
every element in $\de a_1 \w \ldots \w \de a_r\w \de T$ which is rather easy to see and left to the reader.
\be
\end{prf}

The description of $\ann\Omega^n_F(\de S_1 \w \ldots \w \de S_r)$ for 
non-empty sets $S_i \subset F$ is now complete for the cases 
\begin{itemize}
\item $\pdeg_F(S_1 \cup \ldots \cup S_r)=\sum_{i=1}^r \pdeg_F(S_i)$
\item $S_1=\ldots=S_r$
\end{itemize}
which are both extreme cases since for $j \in\{1,\ldots,r\}$ we obviously have 
$$\pdeg_F(S_j) \leq \pdeg_F(S_1 \cup \ldots \cup S_r) \leq \sum_{i=1}^r \pdeg_F(S_i).$$
Proposition \ref{36} gives an insight of the structure of annihilators if the 
$F^p(S_i)$ might have a non-empty intersection. 
Boundaries for the annihilator can be given in all possible cases.

\begin{pro}\label{37}
Let $S_1,\ldots,S_r \subset F$ be non-empty sets with $\pdeg(S_i)=k_i\in\NN$. 
Assume $[F^p(S_1,\ldots,S_i):F^p(S_1,\ldots,S_{i-1})]=p^{\ell_i}$ with $0\leq \ell_i \leq k_i$ for $i=1,\ldots,r$ and let $\{c_{i1},\ldots,c_{i\ell_i}\}\subset F$ be 
$p$-independent such that $F^p(S_1,\ldots,S_i)=F^p(S_1,\ldots,S_{i-1})( c_{i1},\ldots,c_{i\ell_i})$. Now extend these to a $p$-basis 
$\{c_{i1},\ldots,c_{ik_i}\}$ of $F^p(S_i)$. Then we have 
	\begin{align*}
	\sum_{i=1}^r \de c_{i1}\w\ldots\w\de c_{ik_i} \w \Omega^{n-k_i}(F) &\subset \ann\Omega^n_F(\de S_1\w\ldots\w\de S_r) \\
	&\subset \sum_{i=1}^r \de c_{i1}\w\ldots\w\de c_{i\ell_i} \w \Omega^{n-\ell_i}(F).
	\end{align*}
\end{pro}

\begin{prf}
The second inclusion can be shown similarly as in Proposition \ref{36} and the first 
inclusion can be proven in an analogue way as it was done in Proposition \ref{34}.
\be
\end{prf}

Note that the above bounds coincide under the assumptions of Proposition \ref{34} and the 
upper bound becomes $\Omega^n(F)$ if there is $i\in\{1,\ldots,r\}$ with $\ell_i=0$. 
However the simple example below shows that both bounds can easily differ from the actual 
annihilator. The example also shows that
 the condition ``$\pdeg S_i=1$ for $i=1,\ldots,r$'' in Proposition \ref{36} 
is actually necessary for the statement to be true.

\begin{bsp}\label{38}
Set $F:=\FF_p(a,b,c)$ with transcendental elements $a,b,c$. Then it is well known that 
$\{a,b,c\}$ (with ordering $a<b<c$) is a $p$-basis of $F$. Say $r=2$ and $S_1=\{a,b\}, S_2=\{a,c\}$. 
Thus we have $\de S_1\w \de S_2 = \{0,\de a \w \de c, \de b \w \de a, \de b \w \de c\}$, 
hence 
$$\spann_F(\de S_1 \w \de S_2)= \spann_F (\de T \w \de T) \text{ with } T=\{a,b,c\}.$$
So by Lemma \ref{32}(d) the $\Omega$-annihilators of $\de S_1\w \de S_2$ 
and $\de T \w \de T$ coincide and by Proposition \ref{35} we get 
\begin{align*}
\ann\Omega^n_F(\de T \w \de T)= \de a \w \de b \w \Omega^{n-2}(F) +
\de a \w \de c \w \Omega^{n-2}(F) + \de b \w \de c \w \Omega^{n-2}(F).
\end{align*} 
Constructing the bounds of $\ann\Omega^n_F(\de S_1 \w \de S_2)$ as described in 
Proposition \ref{37}, as a lower bound we get
\begin{align*}
\de a \w \de b \w \Omega^{n-2}(F) + \de a \w \de c \w \Omega^{n-2}(F)
\subset \ann\Omega^n_F(\de S_1 \w \de S_2) 
%\subset\de a \w \de b \w \Omega^{n-2}(F) + \de c\w \Omega^{n-1}(F)
\end{align*}
which inclusion is strict since for $n=2$ we have $\de b \w \de c \not\in F\de a \w \de b + F \de a \w \de c$. The upper bound from 
Proposition \ref{37} in this example is given by 
$$\de a \w \de b \w \Omega^{n-2}(F) + \de c\w \Omega^{n-1}(F)$$
which is strictly larger than $\ann\Omega^n_F(\de S_1 \w \de S_2)$ for $n=1$.
\end{bsp}

\end{section}

\begin{section}{Annihilators in $\nu_n(F)$}

In this section our main goal is to use the results from Section \ref{sec2} to find the 
corresponding $\nu$-annihilators of the sets discussed above. For that, 
Kato's lemma will be the key engine for most of the upcoming computations, 
but its usage will also lead to the necessary assumption $F^{p-1}=F$. This restriction 
however is not as strict as it may appear, since the case $p=2$ is of special interest 
for us due to the connection between differential forms and bilinear forms over 
the field $F$ (see \cite[Theorem p. 494]{q9}).%, and $F^{p-1}=F$ is obviously true for $p=2$.

To make it easier for the reader we will repeat Kato's Lemma here.

\begin{pro}\label{23}{\cite[Lemma 2]{q4}}
Let $F=F^{p-1}$ and $\pB=\{b_i \mid i\in I\}$ be a well ordered $p$-basis of $F$. For 
all $x \in F^*$ and $\delta \in \Sigma_n$ 
with $(x^p-x) \frac{\de b_{\delta}}{b_{\delta}} \in \de \Omega^{n-1}(F)+\Omega^n_{<\delta}(F)$, there exist $a_j \in F_{\delta(j)}^*$ for $j=1,\ldots,n$ and $v \in \Omega^n_{<\delta}(F)$ with 
$$ x \frac{\de b_{\delta}}{b_{\delta}} =\frac{\de a_1}{a_1}\w\ldots \w \frac{\de a_n}{a_n} + v.$$
\end{pro}

For the computation of the $\nu$-annihilators, we will need to intersect the 
found $\Omega$-annihilators from Section \ref{sec2} with the group $\nu_n(F)$, which will be done using the following 
two lemmas. We will use $[ \,\cdot\, ]$ as a notation for additive generating systems. Recall that 
$\nu_n(F)= \left[ \frac{\de x_1}{x_1} \w \ldots \w \frac{\de x_n}{x_n} \mid x_i \in F^* \right]$.

\begin{lem}\label{41}
Let $\{a_1,\ldots,a_r,e_1,\ldots,e_{\ell}\}\subset F$ be $p$-independent with $r,\ell\in\NN_0$ 
and assume $F^{p-1}=F$. Then for $t \in \{0,\ldots,\ell\}$
\begin{align*}
	&\ \ \ \left(\sum_{i=1}^r \de a_i \w \Omega^{n-1}(F) + 	
	\sum_{\substack{\{i_1,\ldots,i_t\}\subset\{1,\ldots,\ell\} \\ i_1<\ldots <i_t  }} \de e_{i_1} \w \ldots \w\de e_{i_t}\w \Omega^{n-t}(F) 
	\right) \cap \nu_n(F) 
	\\
	&=\left[ \frac{\de x}{x} ~\Big\vert~ x\in F^p(a_1,\ldots,a_r)^* \right]\w \nu_{n-1}(F)  \\
	&\quad+\ \left[ \frac{\de y_{1}}{y_1}\w\ldots\w\frac{\de y_{t}}{y_{t}} ~\Big\vert~ y_1,\ldots,y_{t} \in F^p(a_1,\ldots,a_r,e_1,\ldots,e_{\ell})^*
	 \right]\w \nu_{n-t}(F).
	\end{align*} 
	where we define $\de e_{i_1} \w \ldots \w\de e_{i_t} =0$ for $t =0$. \\	
In particular for $\ell=0$ (and thus $t=0$) we have 
	\begin{align*}
	\left(\sum_{i=1}^r \de a_i \w \Omega^{n-1}(F)  \right) \cap \nu_n(F) 
	=\left[ \frac{\de x}{x} ~\Big\vert~ x\in F^p(a_1,\ldots,a_r)^* \right]\w \nu_{n-1}(F) 
	\end{align*} 
and for $r=0$ we have 
	\begin{align*}
	\!&\left(\sum_{\substack{\{i_1,\ldots,i_t\}\subset\{1,\ldots,\ell\} \\ i_1<\ldots <i_t  }} \de e_{i_1} \w \ldots \w\de e_{i_t}\w \Omega^{n-t}(F) 
	\right) \cap \nu_n(F) \\
	&=\left[ \frac{\de y_{1}}{y_1}\w\ldots\w\frac{\de y_{t}}{y_{t}} \;\Big\vert\; y_1,\ldots,y_{t} \in F^p(e_1,\ldots,e_{\ell})^*\right]\w \nu_{n-t}(F).
	\end{align*} 	
\end{lem}

\begin{prf}
For the inclusion $(\supseteq)$ let $v=\frac{\de x}{x}\w \chi \in \nu_n(F)$ be a generator with 
$x\in F^p(a_1,\ldots,a_r)^*$ and $\chi \in \nu_{n-1}(F)$. Since $\{a_1,\ldots,a_r\}$ 
is a $p$-basis of $F^p(a_1,\ldots,a_r)$, we find (unique) $\lambda_1,\ldots,\lambda_r \in F$ with $\de x= \sum_{i=1}^r \lambda_i \de a_i$. Thus we have
$$v=\left( \sum_{i=1}^r \frac{\lambda_i}{x} \de a_i \right) \w \chi = \sum_{i=1}^r  \de a_i \w \left( \frac{\lambda_i}{x} \chi\right) 
	\in  \left(\sum_{i=1}^r \de a_i \w \Omega^{n-1}(F)  \right) \cap \nu_n(F).$$
Now let $v=\frac{\de y_{1}}{y_1}\w\ldots\w\frac{\de y_{t}}{y_{t}} \w \eta \in \nu_n(F)$ be a generator of the second summand 
with $y_j \in F^p(a_1,\ldots,a_r,e_1,\ldots,e_{\ell})^*$ 
and $\eta \in \nu_{n-t}(F)$. As before we find representations
$$\de y_j = \sum_{i=1}^r \lambda_{ji} \de a_i + \sum_{m=1}^{\ell} \mu_{jm} \de e_m \ \ \text{ with } \lambda_{ji},\mu_{jm}\in F.$$
Thus we have 
\begin{align*}
	\de y_1\w\ldots\w\de y_{t} =  z + \sum_{j_1=1}^{\ell}\ldots \sum_{j_t=1}^{\ell} \mu_{j_1}\ldots\mu_{j_t} \de e_{j_1}\w\ldots\w\de e_{j_t} 
	\end{align*}
for some $z\in\sum_{i=1}^r \de a_i \w \Omega^{t-1}(F)$, which  
leads to
\begin{align*}
	v&
	= z\w \left((y_1\ldots y_{t})^{-1}\eta\right)  +  \sum_{\substack{\{i_1,\ldots,i_t\}\subset\{1,\ldots,\ell\} \\ i_1<\ldots <i_t  }} \de e_{i_1} \w \ldots \w\de e_{i_t} \w \left(\mu_{j_1}\ldots\mu_{j_t}(y_1\ldots y_{t})^{-1}\eta \right)  
	\\
	&\in \left(\sum_{i=1}^r \de a_i \w \Omega^{n-1}(F) + 	
	\sum_{\substack{\{i_1,\ldots,i_t\}\subset\{1,\ldots,\ell\} \\ i_1<\ldots <i_t  }} \de e_{i_1} \w \ldots \w\de e_{i_t}\w \Omega^{n-t}(F) 
	\right) \cap \nu_n(F).
	\end{align*}

For the inclusion $(\subseteq)$ extend $\{a_1,\ldots, a_r,e_1,\ldots,e_{\ell}\}$ to 
a $p$-basis $\pB=\{b_i \mid i\in I\}$ of $F$ with an ordering such that $a_1,\ldots, a_r,e_1,\ldots,e_{\ell}$ are the smallest elements in $\pB$ with 
$a_1< \ldots< a_r<e_1<\ldots < e_{\ell}$. To simplify notations set 
$$ \mathfrak A:= \left(\sum_{i=1}^r \de a_i \w \Omega^{n-1}(F) + 	
	\sum_{\substack{\{i_1,\ldots,i_t\}\subset\{1,\ldots,\ell\} \\ i_1<\ldots <i_t  }} \de e_{i_1} \w \ldots \de e_{i_t}\w \Omega^{n-t}(F) 
	\right) \cap \nu_n(F)$$
and note that $\mathfrak A$ is an additive group.
%because both $\nu_n(F)$ and $\sum_{i=1}^r \de a_i \w \Omega^{n-1}(F) + \de e_1 \w \ldots \w \de e_{\ell}\w \Omega^{n-\ell}(F)$ are additive groups. 
Now let $\omega \in \mathfrak A$. Since we only work with one element $\omega$ and apply a 
finite number of computation 
steps to this one element, without loss of generality we may 
assume $|\pB|<\infty$ which also implies $|\Sigma_n|< \infty$. Representing $\omega$ in the basis $\bigwedge^n_{\pB}$ and 
assuming $\max_{\pB}(\omega)=\delta \in \Sigma_n$, we get 
\begin{align}\label{e1}
	\omega= \sum_{\sigma\leq \delta} x_{\sigma} \frac{\de b_{\sigma}}{b_{\sigma}}= x_{\delta}\frac{\de b_{\delta}}{b_{\delta}} + \omega_{<\delta}\ \ 
	\text{ with } x_{\sigma}\in F \text{ and } \omega_{<\delta}\in\Omega^n_{<\delta}(F).
	\end{align}
Since $\wp$ respects the filtration on $\Omega^n(F)$ and $\wp(\omega)\in\de \Omega^{n-1}(F)$ by the choice of $\omega$, we have
\begin{align*}
\wp(x_{\delta})\frac{\de b_{\delta}}{b_{\delta}} \in \de \Omega^{n-1}(F) + \Omega^n_{<\delta}(F)
\end{align*}
which allows us to use Kato's lemma \ref{23} so that we get
	\begin{align}\label{e2}
	x_{\delta}\frac{\de b_{\delta}}{b_{\delta}}
	= \frac{\de y_1}{y_1}\w\ldots\w\frac{\de y_n}{y_n} + v_{<\delta} \ \ \text{ with } y_i\in F_{\delta(i)}^* \text{ and } v_{<\delta}\in\Omega^n_{<\delta}(F).
	\end{align}
Now due to the assumptions, at least one of the following two cases must 
occur.
\begin{enumerate}[(i)]
\item There is $ i \in \{1,\ldots,r\}$ with $\de a_i \mid \de b_{\delta}$
\item There is a subset $\{i_1,\ldots,i_t\} \subset \{1,\ldots, \ell\}$ with $i_1 < \ldots < i_t$ such that $ \de e_{i_1} \w \ldots \w \de e_{i_t} \mid \de b_{\delta}$
\end{enumerate}
First assume case (i) and say $\de a_s \mid \de b_{\delta}$ with $s \in \{1,\ldots,r\}$ 
and further assume $\de a_1, \ldots , \de a_{s-1} \nmid \de b_{\delta}$. 
Using the ordering on $\pB$, we have $\delta(1)=s$ and $y_1 \in F^p(a_1,\ldots,a_s)^* 
 \subset F^p(a_1,\ldots,a_r)^*$. If we now insert Equation \eqref{e2} into \eqref{e1}, 
 we get 
	\begin{align*}
	\omega= \frac{\de y_1}{ y_1}\w\ldots\w\frac{\de y_n}{y_n} +\omega_{<\delta} + v_{<\delta}
	\end{align*}
and since we have $\frac{\de y_1}{ y_1}\w\ldots\w\frac{\de y_n}{y_n}\in \mathfrak A$ by the 
already shown first inclusion, we also 
get $\omega_{<\delta} + v_{<\delta} = \omega - \frac{\de y_1}{ y_1}\w\ldots\w\frac{\de y_n}{y_n} \in \mathfrak A$ and we
may proceed by induction on $\max_{\pB}(\omega)$.
	
Now assume case (ii), i.e. say $\de e_{i_1} \w \ldots \w \de e_{i_t} \mid \de b_{\delta}$. 
%and further assume that $(i_1,\ldots,i_t)$ is lexicographically the smallest element 
%in $\{1, \ldots, \ell\}^t$ fulfilling this property. 
Due to case (i) already being treated, we may additionally assume that none of the forms 
$\de a_1,\ldots, \de a_r$ divide $\de b_{\delta}$. Thus, by the ordering of $\pB$  and from \eqref{e2} we get
\begin{align*}
	y_1,\ldots,y_{t}\in F_{\delta(i_t)} \subset F^p(a_1,\ldots,a_r,e_1,\ldots e_{\ell})^*.
	\end{align*}
Combining Equations \eqref{e2} and \eqref{e1} once more, and repeating the same arguments 
as in case (i), we may conclude the proof by induction on $\max_{\pB}(\omega)$.	
\be
\end{prf}

Note that in Lemma \ref{41}, the $y_i$ might also be chosen in $F^p(a_1,\ldots,a_r,e_1,\ldots,e_{\ell})\setminus F^p(a_1,\ldots,a_r)$ because if we have $y_i \in  F^p(a_1,\ldots,a_r)$, the resulting generator $\frac{\de y_1}{y_1}\w\ldots\w\frac{\de y_{t}}{y_{t}}$ (after possibly reordering the slots) can be viewed as an element of $\left[ \frac{\de x}{x} \mid x\in F^p(a_1,\ldots,a_r)^* \right]\w \nu_{n-1}(F)$.

If we now combine the results from Section \ref{sec2} with Lemma \ref{41}, 
we are able to give a precise description for the $\nu$-annihilators of the sets 
studied in Section \ref{sec2}.

\begin{theo}\label{43}
Assume $F^{p-1}=F$.
	\begin{enumerate}[(a)]
	\item Let $S_1,\ldots,S_{r+1}\subset F$ be non-empty sets with $\pdeg S_i=1$ for $i=1,\ldots,r$ and $\de S_1 \w \ldots \de S_r \w \de S_{r+1}\neq \{0\}$. 
Set $\pdeg(S_1 \cup \ldots \cup S_r \cup S_{r+1})=r+\ell$ with $\ell \in \ZZ$. Then $\ell \geq 1$ and for any $p$-independent 
$\{a_1,\ldots,a_r,e_1,\ldots, e_{\ell}\}\subset F$ with $F^p(S_i)=F^p(a_i)$ for $i=1,\ldots,r$ 
and $F^p(a_1,\ldots, a_r)(S_{r+1})=F^p(a_1,\ldots,a_r)(e_1,\ldots,e_{\ell})$ we have 
		\begin{align*}
		&\ann\nu^n_F(\de S_1 \w \ldots \w \de S_r \w \de S_{r+1}) =\bigg[ \frac{\de x}{x}  ~\Big\vert~ x\in F^p(a_1,\ldots,a_r)^* \bigg]\w \nu_{n-1}(F) \\
		&\ \ + \bigg[\frac{\de y_1}{y_1}\w\ldots\w\frac{\de y_{\ell}}{y_{\ell}} ~\Big\vert~ 
		y_1,\ldots,y_{\ell}\in F^p(a_1,\ldots,a_r,e_1,\ldots e_{\ell}) \bigg]\w\nu_{n-\ell}(F).
		\end{align*}
	\item Let $S\subset F$ be a non-empty set with $\pdeg (S)=k\in\NN$ and let $\{a_1,\ldots,a_k\}$ be a $p$-basis of $F^p(S)$. 
For $r\in\{1,\ldots,k\}$ set $t=k-r+1$. Then 
		\begin{align*}
		\ann\nu^n_F\left({\bigwedge}^r\de S \right) = \bigg[ \frac{\de y_1}{y_1}\w\ldots\w\frac{\de y_t}{y_t} ~\Big\vert~ y_1,\ldots,y_t\in F^p(a_1,\ldots,a_k)^*  \bigg] \w\nu_{n-t}(F),
		\end{align*}
	and for $r>k$ we have $\ann\nu^n_F({\bigwedge}^r\de S)=\nu_n(F)$.
\end{enumerate}
\end{theo}

\begin{rem}
\begin{enumerate}[(a)]
	\item Unfortunately our methods do not apply for the computation of the intersection 
	of $\sum_{i=1}^r \de a_{i1}\w\ldots\w\de a_{ik_i} \w \Omega^{n-k_i}(F)$ with 
	$\nu_n(F)$, even if the set $\{a_{ij} \mid i=1,\ldots,r\ j=1,\ldots, k_i\}$ is 
	$p$-independent. However a general additive behavior of the $\nu$-annihilator 
	as described in \ref{34} for the $\Omega$-annihilator can be easily disproved 
	by a counterexample.

	\item In Theorem \ref{43}(a), if we additionally assume $\pdeg (S_{r+1})=1$ and 
	rename $e_1$ as $a_{r+1}$, we get 
	\begin{align*}
	\ann\nu^n_F(\de S_1 \w \ldots \w \de S_r \w \de S_{r+1})=\bigg[ \frac{\de x}{x}  ~\Big\vert~ x\in F^p(a_1,\ldots,a_{r+1})^* \bigg]\w \nu_{n-1}(F).
	\end{align*}
	Under these assumptions and additionally assuming $p=2$ and $S_i=\{a_i\}$, this annihilator was already known and can be found in \cite[Theorem 4.1]{q8} but was derived 
	using quite different methods.
	
	\item The case $r=1$ in Theorem \ref{43}(b) gives us 
	 	\begin{align*}
	\ann\nu^n_F(\de S)= \bigg[ \frac{\de y_1}{y_1}\w\ldots\w\frac{\de y_k}{y_k} ~\Big\vert~ y_1,\ldots,y_k\in F^p(a_1,\ldots,a_k)^*  \bigg] \w\nu_{n-k}(F).
	\end{align*}
	This annihilator, even if not stated directly, can be easily derived from the results 
	given in \cite{q2} (combine Proposition 7.2 with the 
	arguments of Section 9).

\end{enumerate}

\end{rem}
\end{section}

\begin{section}{Annihilators and kernels of algebraic field extensions}

Apart from being of interest of its own, annihilators seem to have an interesting connection to specific kernels of field extensions 
and can be used to determine yet unknown kernels. In this last section, our goal is to present this connection and interpret some 
known kernels as annihilators. 

To start of, for a field extension $K/F$ we will shortly write $\omega_K$ for a form $\omega\in \Omega^n(F)$ interpreted as a form in $\Omega^n(K)$ 
and define 
\begin{align*}
	\Omega^n(K/F)&=\{ \omega \in \Omega^n(F) \mid \omega_K=0 \text{ in } \Omega^n(K)   \}\\
	\nu_n(K/F)&=\{ \chi \in \nu_n(F) \mid \chi_K=0 \text{ in } \nu_n(K)   \}= \Omega^n(K/F) \cap \nu_n(F).
\end{align*}
Furthermore let $\overline F$ be an algebraic closure of $F$ and $\alpha \in \overline F$. We write  $\min_{F,\alpha}$
for the monic irreducible polynomial with root $\alpha$ over $F$. Set
\begin{align*}
\mathcal C_F(\alpha):= \Big\{ c_i \in F \mid \min\nolimits_{F,\alpha}=\sum_{i=0}^k c_iX^i  \Big\} \setminus \{0\}.
\end{align*}

With this the kernel of an arbitrary simple algebraic field extension can be described as follows.

\begin{pro}\label{51}{\cite[Proposition 7.2]{q2}}
Let $\overline F$ be an algebraic closure of $F$, $\alpha \in \overline F$ and set $K=F(\alpha)$.
\begin{enumerate}[(a)]
	\item If $\alpha$ is separable over $F$, then $\Omega^n(K/F)=\{0\}=\nu_n(K/F)$.
	\item If $\alpha$ is not separable over $F$, then 
	$$\Omega^n(K/F)=\ann\Omega^n_F(\de \mathcal C_F(\alpha) )\ \ \text{ and } \ \ \nu_n(K/F)=\ann\nu^n_F(\de \mathcal C_F(\alpha) )$$
	and both annihilators can be precisely described using Proposition \ref{35} and Theorem \ref{43}.
\end{enumerate}
\end{pro}

Since the case of a simple algebraic extension is covered with Proposition \ref{51}, let us now move on to 
non-simple extensions. For that recall that if $\{ b_1,\ldots, b_r \}\subset F$ is $p$-independent over $F$, 
then for any $m_1,\ldots,m_r\in \NN$ the set $\{ \sqrt[p^{m_1}]{b_1}, \ldots, \sqrt[p^{m_r}]{b_r}\}$ is $p$-independent over 
 $F\left(\sqrt[p^{m_1}]{b_1}, \ldots, \sqrt[p^{m_r}]{b_r} \right)$. With this it is an easy exercise to prove the following well known lemma.

\begin{lem}\label{52}
Let $\{ b_1,\ldots, b_r \}\subset F$ be $p$-independent and $m_1,\ldots,m_r\in\NN$. Let $\pB$ be a $p$-basis of $F$ containing the $b_i$. 
Set $E=F\left(\sqrt[p^{m_1}]{b_1}, \ldots, \sqrt[p^{m_r}]{b_r} \right)$, $\beta_i=\sqrt[p^{m_i}]{b_i}$, $T=\bigtimes_{i=1}^r\{1,\ldots,p^{m_i}-1\}$ 
and $\mathcal B'=\pB \setminus \{ b_1,\ldots, b_r \}$. Then $\Omega^n(E)$ is given by 
\begin{align*}
\Omega^n(E)= \bigoplus_{\substack{(i_1,\ldots,i_r)\in T\\ \{\beta_{j_1},\ldots,\beta_{j_s}\}\subset\{\beta_1,\ldots,\beta_r\}\\ \text{for }s=0,\ldots,r \text{ with } j_1<\ldots <j_s}}
	\beta_1^{i_1}\ldots \beta_r^{i_r} \de \beta_{j_1}\w\ldots \w \de \beta_{j_s}\w \Omega^{n-s}_{\pB'}(F)
\end{align*}
where the direct sum is taken as a sum of $F$-vector spaces and both sides of the equation are viewed as $F$-vector spaces.
\end{lem}

So with this, we will now take a closer look on purely inseparable extensions and start with so called modular purely inseparable 
extensions, i.e.  with purely inseparable extensions $E/F$, such that we can find $p$-independent 
elements $b_1,\ldots,b_r \in F$ and $m_1,\ldots,m_r\in \NN$ with $E=F\left(\sqrt[p^{m_1}]{b_1}, \ldots, \sqrt[p^{m_r}]{b_r} \right)$. 
Using \cite[Theorem 4.1]{q55} and combining it with Proposition \ref{34} resp. with Theorem \ref{43}, we get the following result.

\begin{pro}\label{53}
Let $\{ b_1,\ldots, b_r \}\subset F$ be $p$-independent, $m_1,\ldots,m_r\in\NN$ and set $E=F\left(\sqrt[p^{m_1}]{b_1}, \ldots, \sqrt[p^{m_r}]{b_r} \right)$. 
Then 
\begin{align*}
\Omega^n(E/F)&= \sum_{i=1}^r \de b_i \w \Omega^{n-1}(F)= \ann\Omega^n_F(\de b_1 \w \ldots \w \de b_r), \\
\nu_n(E/F)&=\ann\nu^n_F(\de b_1 \w \ldots \w \de b_r)
\end{align*}
and for $F^{p-1}=F$ we have
\begin{align*}
\nu_n(E/F)= \left[\frac{\de x}{x} \mid x \in F^p(b_1,\ldots,b_r)^* \right]\w \nu_{n-1}(F).
\end{align*}
\end{pro}

For non-modular purely inseparable extensions, the kernels are far more difficult to determine. However with 
the use of annihilators, we can get a first insight on the structure of these kernels.

\begin{theo}\label{54}
Let $b_1,\ldots,b_r,b\in F\setminus F^p$ be such that $\{b_1,\ldots,b_r\}$ is $p$-independent and $b\in F^p(b_1,\ldots,b_r)$. For 
$m_1,\ldots,m_r,m \in \NN$ with $m \leq m_1,\ldots,m_r$ set $E=F(\sqrt[p^{m_1}]{b_1}, \ldots, \sqrt[p^{m_r}]{b_r}, \sqrt[p^{m}]{b} )$.
\begin{enumerate}[(a)]
	\item If there exists $t\in\NN$, $t\geq m$ with $b \in F^{p^t}(b_1,\ldots,b_r)$, then $E/F$ is modular 
	with $E=F\left(\sqrt[p^{m_1}]{b_1}, \ldots, \sqrt[p^{m_r}]{b_r} \right)$ 
	and the kernels $\Omega^n(E/F)$ and $\nu_n(E/F)$ can be described using Proposition \ref{53}.
	\item Assume (a) does not hold and choose $t\in\{1,\ldots,m-1\}$ maximal with $b\in F^{p^t}(b_1,\ldots,b_r)$. 
	Write $b=\sum_{i=(i_1,\ldots,i_r)\in T} x_i^{p^t} b_1^{i_1}\ldots b_r^{i_r}$ with $T=\{0,\ldots,p^t-1\}^r$ and set $S:=\{ x_i \mid i\in T \}$. 
	Then 
	\begin{align*}
	\Omega^n(E/F)&=\ann\Omega^n(\de b_1 \w \ldots \w \de b_r \w \de S) \ \ \text{ and} \\
	\nu_n(E/F)&=\ann\nu^n_F(\de b_1 \w \ldots \w \de b_r \w \de S)
	\end{align*}
and both annihilators can be precisely described using Proposition \ref{36} and Theorem \ref{43}.
\end{enumerate}

\end{theo}

\begin{prf}
Set $M=F\left(\sqrt[p^{m_1}]{b_1}, \ldots, \sqrt[p^{m_r}]{b_r} \right)$. Let us start by proving (a) and assume we find $t \geq m$ with  
$$b=\sum_{i=(i_1,\ldots,i_r)\in T} x_i^{p^t} b_1^{i_1}\ldots b_r^{i_r} \in F^{p^t}(b_1,\ldots,b_r).$$
Then by $m\leq m_1,\ldots,m_r$ and $t-m > 0$ we get 
	\begin{align*}
	\sqrt[p^m]{b}= 
	\sum_{(i_1,\ldots,i_r)\in  T} x_{i}^{p^{t-m}} \left(\sqrt[p^{m}]{b_1}\right)^{i_1}\ldots\left(\sqrt[p^{m}]{b_r}\right)^{i_r} 
	\in M
	\end{align*}
from which we conclude $E=M$ as claimed. 

Let us now consider the case (b) and say we have a maximal $t\in\{1,\ldots,m-1\}$ with $b\in F^{p^t}(b_1,\ldots,b_r)$. We start by showing the  
inclusion $(\supseteq)$. 
Assume $\de b_1 \w \ldots \w \de b_r \w \de S =\{0\}$. By the $p$-independence of the $b_1,\ldots,b_r$, we get 
$S \subset F^p(b_1,\ldots,b_r)$, so for every $x_i \in S$ we find $z_{ij}\in F$, $j = (j_1,\ldots,j_r)\in \{0,\ldots,p-1\}^r$ with 
$$ x_i= \sum_{j=(j_1,\ldots,j_r)\in \{0,\ldots,p-1\}^r} z_{ij}^p b_1^{j_1}\ldots b_r^{j_r}.$$
Inserting this in the representation of $b$ we get 
	\begin{align*}
	b=\sum_{i\in T} x_{i}^{p^t} b_1^{i_1}\ldots b_r^{i_r} 
	= \sum_{i\in T} \sum_{j\in \{0,\ldots,p-1\}^r} z_{ij}^{p^{t+1}}    b_1^{i_1+j_1}\ldots b_r^{i_r+j_r} \in F^{p^{t+1}}(b_1,\ldots,b_r)
	\end{align*}
which contradicts $t$ being maximal. Thus we have $\de b_1 \w \ldots \w \de b_r \w \de S \neq \{0\}$. 

Now choose $e_1,\ldots, e_{\ell}\in F$ such that $\{b_1,\ldots,b_r,e_1,\ldots,e_{\ell}\}$ is $p$-independent with 
$F^p(b_1,\ldots,b_r)(S)=F^p(b_1,\ldots,b_r)(e_1,\ldots,e_{\ell})$. By Proposition \ref{36}, to prove the first inclusion it is enough to 
check that all of the forms $\de b_1,\ldots, \de b_r$ and $\de e_1 \w \ldots \w \de e_{\ell}$ become zero forms over $E$. 
Since $b_1, \ldots, b_r \in E^p$, this is obvious for the forms $\de b_1,\ldots, \de b_r$. So let us now check the last form. Because of 
$t< m \leq m_1,\ldots,m_r$, we get $\sqrt[p^{t}]{b}, \sqrt[p^{t}]{b_1}, \ldots, \sqrt[p^{t}]{b_r} \in E^p$ and over $E$ we have 
\begin{align}\label{e6}
	0&=\de\left( \sqrt[p^t]{b}\right)  = \de \left(\sum_{(i_1,\ldots,i_r)\in T} x_{i} \left(\sqrt[p^t]{b_1}\right)^{i_1}\ldots \left(\sqrt[p^t]{b_r}\right)^{i_r} \right) \\ 
	&= \sum_{(i_1,\ldots,i_r)\in T}  \left(\sqrt[p^t]{b_1}\right)^{i_1}\ldots \left(\sqrt[p^t]{b_r}\right)^{i_r} \de x_{i}. \nonumber
	\end{align}
Since $x_i \in S \subset S\cup \{b_1,\ldots,b_r\}\subset F^p(b_1,\ldots,b_r,e_1,\ldots,e_r)$, we find (unique) $\lambda_{ik},\mu_{iq}\in F$ for $k=1,\ldots,r$, $q=1,\ldots,\ell$ with 
\begin{align}\label{e7}
\de x_i = \sum_{k=1}^r \lambda_{ik}\de b_k + \sum_{q=1}^{\ell} \mu_{iq} \de e_q.
\end{align}
Inserting this representation of $\de x_i$ into Equation \eqref{e6} and using $(\de b_k)_E=0$, we get 
	\begin{align}\label{e8}
	0= \sum_{q=1}^{\ell} \left( \sum_{(i_1,\ldots,i_r)\in T}\left(\sqrt[p^t]{b_1}\right)^{i_1}\ldots \left(\sqrt[p^t]{b_r}\right)^{i_r}\mu_{iq} \right) \de e_q
	\end{align}
over $E$. Now assume $\mu_{iq}=0$ for all $i \in T$ and $q \in \{1,\ldots,\ell\}$. Thus by Equation \eqref{e7} we 
get $\de x_i \in \spann_F(\de b_1,\ldots,\de b_r)$ hence $x_i \in F^p(b_1,\ldots,b_r)$ for all $i\in T$ but this was already shown to be wrong at the beginning 
of the proof. So there is at least one $\mu_{iq}\neq 0$. Since $\{b_1,\ldots,b_r\}$ is $p$-independent, the elements  
$\left(\sqrt[p^t]{b_1}\right)^{i_1}\ldots \left(\sqrt[p^t]{b_r}\right)^{i_r}$ with $(i_1,\ldots,i_r)\in T$ are all pairwise distinct and part of the canonical $F$-basis of 
$M$, hence $F$-linear independent. Thus there is at least one $q \in \{1,\ldots,\ell\}$ with $\sum_{(i_1,\ldots,i_r)\in T}\left(\sqrt[p^t]{b_1}\right)^{i_1}\ldots \left(\sqrt[p^t]{b_r}\right)^{i_r}\mu_{iq} \neq 0$ and Equation \eqref{e8} gives a non-trivial linear representation of zero 
using the vectors $\de e_1,\ldots,\de e_{\ell}$ over $E$ which means the forms $\de e_1,\ldots, \de e_{\ell}$ are linearly dependent over $E$ and by 
Lemma \ref{21} we finally see $(\de e_1 \w \ldots \w \de e_{\ell})_E=0$.

Let us now prove the inclusion $(\subseteq)$. Set $U=\bigtimes_{j=1}^r \{0,\ldots,p^{m_j}-1\}$. We obviously have 
$E=M(\sqrt[p^m]{b})=M(\sqrt[p^{m-t}]{b_0})$ with 
\begin{align}\label{e9}
b_0=\sqrt[p^t]{b}=\sum_{(i_1,\ldots,i_r)\in T} x_{i} \big(\sqrt[p^t]{b_1}\big)^{i_1}\ldots \big(\sqrt[p^t]{b_r} \big)^{i_r}\in M.
\end{align}
Assume $b_0 \in M^p$, then for suitable $z_u\in F$, $u=(u_1,\ldots,u_r)\in U$, we have 
\begin{align*}
b_0=\left( \sum_{u\in U} z_u \left(\sqrt[p^{m_1}]{b_1}\right)^{u_1}\ldots \left(\sqrt[p^{m_r}]{b_r}\right)^{u_r}  \right)^p 
	= \sum_{u\in U} z_u^p \left(\sqrt[p^{m_1}]{b_1}\right)^{pu_1}\ldots \left(\sqrt[p^{m_r}]{b_r}\right)^{pu_r}
\end{align*}
and inserting this into the relation $b_0^{p^t}=b$, we get 
	\begin{align}\label{e10}
	b=b_0^{p^t}=\sum_{u\in U} z_u^{p^{t+1}} \left(\sqrt[p^{m_1}]{b_1}\right)^{p^{t+1}u_1}\ldots \left(\sqrt[p^{m_r}]{b_r}\right)^{p^{t+1}u_r}.
	\end{align}
Comparing coefficients in Equation \eqref{e10} then leads to $b\in F^{p^{t+1}}(b_1,\ldots,b_r)$ which contradicts the maximality of $t$. Thus    
$b_0 \in M \setminus M^p$ and in particular we have $[E:F]=p^{m_1+\ldots+ m_r + m -t}$ as well as 
$\Omega^n(E/M)= \ann\Omega^n_M(\de b_0)$ by Proposition \ref{53}. Now let $\omega\in \Omega^n(E/F)$. Then $\omega_M \in \Omega^n(E/M)$ and 
$\omega_M \w \de b_0=0$ in $M$ so by using the representation of $b_0$ given in \eqref{e9} together with 
$t<m_1,\ldots,m_r$, we get 
	\begin{align*}
	0=\omega_M \w \de b_0 
	= \sum_{(i_1,\ldots,i_r)\in T}  \left(\sqrt[p^t]{b_1}\right)^{i_1}\ldots \left(\sqrt[p^t]{b_r}\right)^{i_r}  \omega_M \w \de x_{i}.
	\end{align*}
and Lemma \ref{52} applied to $M$ then leads to 
	\begin{align*}
 \omega_M \w \de x_i =0 \ \ \text{ for all } i=(i_1,\ldots,i_r)\in T.
	\end{align*}
Since $x_i \in F$ for all $i=(i_1,\ldots,i_r)\in T$, we have 
\begin{align*}
\omega \w \de x_i \in \Omega^{n+1}(M/F) \overset{\ref{53}}{=}\ann\Omega^{n+1}_F( \de b_1 \w \ldots \w \de b_r) \ \text{ for all } i=(i_1,\ldots,i_r)\in T
\end{align*}
which finally shows 
\begin{align*}
\omega \w \de x_i \w \de b_1 \w \ldots \w \de b_r=0 \text{ over } F \ \text{ for all } i=(i_1,\ldots,i_r)\in T.
\end{align*}
\be
\end{prf}

The kernel $\Omega^n(E/F)$ stated in Theorem \ref{54} in the case $p=2$ was simultaneously and independently discovered by Aravire, Laghribi and O'Ryan and and can be found in
	\cite[Theorem 5.5]{q56} where the same kernel is given using a different representation 
	without usage of annihilators (see for comparison the PhD thesis of the author from 2017 \cite[Theorem 6.5]{q57}). But it appears that the representation of $\Omega^n(E/F)$ given in \cite[Theorem 5.5]{q56} is not applicable to all cases since the coefficients used in the description of the kernel do not need to lie in the basefield $F$.

\begin{rem}\label{542}
\begin{enumerate}[(a)]
\item Note that the condition $m \leq m_1,\ldots,m_r$ becomes obsolete in Theorem \ref{54} if $\de b = \sum_{j=1}^r \lambda_j \de b_j$ with $\lambda_j \neq 0$ for all $j$. In this 
	case $(\{b_1,\ldots,b_r\}\setminus \{b_j\})\cup \{b\}$ is $p$-independent for all $j\in\{1,\ldots ,r\}$ with $b_j \in F^p(b,b_1, \ldots, b_{j-1},b_{j+1},\ldots, b_r)$. 
	So in this case the condition $m \leq m_1,\ldots,m_r$ may be assumed without loss of generality. Also note that this condition can always be assumed to be true for $r=1$ (see (\ref{b})).

\item\label{b} We would like to point out that Theorem \ref{54} gives a complete characterization of kernels $\Omega^n(E/F)$ for purely inseparable extensions 
generated by two elements. The modular case is already known and for the non-modular case note that for a $p$-dependent set $\{a,b\}\subset F\setminus F^p$ and $E=F\left(\sqrt[p^t]{a},\sqrt[p^s]{b}\right)$ we may assume $s\leq t$ without loss of generality, 
since $a\in F^p(b)$ is equivalent to $b \in F^p(a)$.
\item  Note that the condition $m \leq m_1,\ldots,m_r$ in Theorem \ref{54} is actually necessary for the statement to be true. To see that, assume that 
	Theorem \ref{54} is true without the restriction $m \leq m_1,\ldots,m_r$. Now take the field 
	$F=\FF_2(a,b,c)$ of characteristic $p=2$ for some transcendental elements $a,b,c$. We further set $E=F\left( \sqrt[4]{a}, \sqrt{b}, \sqrt[8]{c^4a}\right)$. 
	Obviously $\{a,b\}$ is $2$-independent with $c^4a\in F^4(a,b)\setminus F^8(a,b)$. So case (b) of Theorem \ref{54} would apply and give us
	$$\Omega^n(E/F)=\ann\Omega^n(\de a \w \de b \w \de c)=\de a \w\Omega^{n-1}(F)  +\de b \w\Omega^{n-1}(F)+\de c \w\Omega^{n-1}(F).$$
	But we also have $E=F\left( \sqrt[8]{c^4a}, \sqrt{b},\sqrt[4]{a} \right)$ with $\{c^4a,b\}$ being $2$-independent and $a\in F^4(c^4a,b)$. So now case (a) would apply 
	and we get 
	$$\Omega^n(E/F)=\ann\Omega^n(\de \;(c^4a) \w \de b)=\de a \w\Omega^{n-1}(F)  +\de b \w\Omega^{n-1}(F)$$
	which differs from the first description of $\Omega^n(E/F)$ above for $n=1$.

\end{enumerate}

\end{rem}

\begin{rem}\label{55}
\begin{enumerate}[(a)]	
	\item Using the results as well as the notation from Proposition \ref{53}, we see that 
	\begin{align*}
	\Omega^n(E/F) = \ann \Omega^n_F\left( \de \mathcal{C}_F \left(\sqrt[p^{m_1}]{b_1}\right) \w \ldots \w\de \mathcal C_F \left(\sqrt[p^{m_r}]{b_r}\right) \right).
	\end{align*}
	Set $\beta_i = \sqrt[p^{m_i}]{b_i}$,  $F_i=F(\beta_1,\ldots, \beta_{i})$ for $i=1,\ldots, r$ and $F_0=F$. 
	Note that since the elements $b_1,\ldots, b_r$ are $p$-independent, we have $\min_{F,\beta_i}=X^{p^{m_i}}-b_i\in F[X]$ and 
	$\min_ {F,\beta_i}=\min_{F_{i-1},\beta_i}$ for $i=1,\ldots,r$ (To be precise, this property is equivalent for the $b_1,\ldots,b_r$ to be $p$-independent). 
	Thus we can also rewrite the kernel as  
	\begin{align*}
	\Omega^n(E/F) = \ann \Omega^n_F\left( \de \mathcal{C}_{F_0} \left(\beta_1\right) \w \ldots  \w\de \mathcal C_{F_{r-1}} \left(\beta_r\right) \right).
	\end{align*}
	\item Now we use the results and the notations from the proof of Theorem \ref{54} and set $\beta_i = \sqrt[p^{m_i}]{b_i}$ for $i=1,\ldots, r$ as well 
	as $\beta=\sqrt[p^{m}]{b}$. Then $\min_{M,\beta}=X^{p^{m-t}}-b_0$, so the set $S$ consists of all the coefficients needed to represent $b_0$ using 
	the canonical $F$-basis $\{ \beta_1^{u_1}\ldots \beta_r^{u_r} \mid (u_1,\ldots,u_r) \in \bigtimes_{i=1}^r \{0,\ldots,p^{m_i}-1\} \}$ of $M$. 
	
	\item Note that if we set $p=2$, then by using Kato's Theorem \cite[Theorem p. 494]{q9} we are able to translate all the results for 
	the $\nu$-kernels and the $\nu$-annihilators to the theory of bilinear forms over $F$. The annihilated sets $\de S_1 \w \ldots \w \de S_r$ with 
	$S_i \subset F$ can be transported as a set of so called bilinear Pfister forms. Since 
	this transfer is a standard procedure by now, we leave it to the reader to fill in the details and refer to \cite[Section 9]{q2} or \cite[Section 5]{q55} for 
	additional comments (or to \cite[Kapitel 8]{q57} for a rather detailed description).
\end{enumerate}
\end{rem}

Comparing the interpretations given in Remark \ref{55}(a)+(b) and reviewing the results above, the following questions come to mind:

\textit{Question 1:} Let $K=F(\alpha_1,\ldots,\alpha_r)/F$ be an algebraic $r$-fold field extension such that $\alpha_i$ is not 
separable over $F$ for $i=1,\ldots,r$ (being an $r$-fold extension means $K$ cannot be written with less than $r$ algebraic generators over $F$). Can we find suitable sets $S_i\subset F$ such that 
$$\Omega^n(K/F)=\ann\Omega^n_F( \de S_1 \w \ldots \w \de S_r)?$$

\textit{Question 2}: Set $K=F(\alpha, \beta)$ such that $K/F$ is 
a $2$-fold extension with $\alpha,\beta$ both being algebraic and 
non-separable over $F$. Set $[F(\alpha):F]=s$ 
and $\mathcal C_{F(\alpha)}(\beta)=\{b_1,\ldots,b_k\}$. Now write $b_j=\sum_{i=1}^s a_{ij}\alpha^{i}$ for $j=1,\ldots,k$ with $a_{ij}\in F$. Do we have 
$$ \Omega^n(K/F)= \ann\Omega^n_F(\de \mathcal C_F(\alpha) \w \de\, \{ a_{ij} \mid i=1,\ldots,s, j=1,\ldots,k\} )?$$ 

%\textit{Question 3:} Can the ideas of Question 2 be generalized to $r$-fold extensions of $F$?

Note that both questions have a positive answer for all kernels which are known so far. 

This paper is based on results of the PhD thesis of the author.

\vspace*{0.5cm}

\noindent\textbf{Acknowledgement}

The author thanks the anonymous referee for various comments and suggestions to improve the presentation of the paper.

\end{section}

%\section*{Acknowledgments and Notes} 
%\nocite{*}
\bibliographystyle{abbrv}
\bibliography{literatur}

\end{document}